%% file: expo-main.tex
\documentclass[11pt]{amsart}

\usepackage{graphicx}
\usepackage{amscd}
\usepackage{amsmath}
\usepackage{amsfonts}
\usepackage{amssymb}
\usepackage{amsxtra}
\usepackage{xspace}
\usepackage{array}

\theoremstyle{plain}

\newtheorem{proposition}{Proposition}

\theoremstyle{definition}

\numberwithin{equation}{section}

\newcommand{\be}{\begin{enumerate}}
\newcommand{\ee}{\end{enumerate}}
\newcommand{\beq}{\begin{equation}}
\newcommand{\eeq}{\end{equation}}
\newcommand{\bprop}{\begin{proposition}}
\newcommand{\eprop}{\end{proposition}}

\newcommand{\rationals}{\mathbb{Q}}
\newcommand{\integers}{\mathbb{Z}}

\newcommand{\A}{\mathcal{A}}

\DeclareMathOperator{\len}{\ell}

\DeclareMathOperator{\htt}{ht} 
\DeclareMathOperator{\wt}{wt}

\newcommand{\tko}[1][]{\mathcal{P}(#1; t)}

\newcommand{\sch}[1][\lambda]{\chi_{#1}}

\newcommand{\kma}{\mathfrak{g}}

\newcommand{\p}[1][\gamma]{\mathrm{Par}(#1)}
\newcommand{\ppi}[1][\gamma]{\mathrm{Par}(#1,\alpha_i)}
\newcommand{\ppni}[1][\gamma]{\mathrm{Par}(#1,\widehat{\alpha_i})}

\newcommand{\ro}{\Delta}
\newcommand{\rp}{\Delta^+}

\newcommand{\highroot}{\tilde{\alpha}} 

\begin{document}

\title[]{A note on exponents vs root heights for complex simple Lie algebras}
\author{Sankaran Viswanath}
\address{Department of Mathematics\\
University of California\\
Davis, CA 95616, USA}
\email{svis@math.ucdavis.edu}
\subjclass[2000]{}

\begin{abstract}
We give an elementary combinatorial proof of a special case of a 
result due to Bazlov and Ion concerning the Fourier coefficients of
the Cherednik kernel. This can be used to 
give yet another  proof of the classical fact that for a
complex
simple Lie algebra $\kma$, the
partition formed by the exponents of $\kma$ is dual to that formed  by the numbers of positive
roots at each height.
\end{abstract}

\maketitle

\input{expo.tex}

\bibliographystyle{amsplain}
\bibliography{expo}
\end{document}

%% file: expo.tex
\section{Introduction}
Let $\kma$ be a finite dimensional, complex simple Lie algebra of rank $n$ with
associated root system $\ro$, simple roots $\alpha_i$ ($i=1\cdots
n$) and set of positive roots $\rp$. Let $Q$ be the root lattice of $\kma$ and $Q^+$ denote the set comprising
$\integers^{\geq  0}$ linear combinations of the $\alpha_i$. For each
$\alpha \in Q$, let $e^\alpha$ denote the corresponding formal
exponential; these satisfy the usual rules: $e^0=1$ and
$e^{\alpha+\beta} = e^{\alpha}\,e^{\beta}$. We define
$\A :=\rationals[t]\,[[e^{-\alpha_1},\cdots,e^{-\alpha_n}]]$. Thus a
typical element of $\A$ is a power series of the form $\sum_{\beta \in Q^+}
c_{\beta}(t) e^{-\beta}$ where each $c_\beta (t) \in \rationals[t]$.

Consider the element $\xi \in \A$ defined by:
\begin{align}
\xi &:= \displaystyle\prod_{\alpha \in \rp} \frac{1-e^{-\alpha}}{1-t e^{-\alpha}}\notag\\
&=  \displaystyle\prod_{\alpha \in \rp} (1 + (t-1)e^{-\alpha} +t(t-1) e^{-2\alpha}+
  t^2(t-1) e^{-3\alpha} + \cdots ) \label{xiq}
\end{align}
Given $\beta = \sum_{i=1}^n b_i \alpha_i \in Q^+$, define its {\em
  height} to be $$ \htt \beta:=\sum_{i=1}^n b_i$$
The main objective of this short note is to give an elementary
combinatorial proof of the following proposition:

\begin{proposition} \label{cl1}
For $\beta \in \rp$, the coefficient of $e^{-\beta}$ in $\xi$ is
$(\,t^{\htt(\beta)}  - t^{\htt(\beta)-1}\,)$
\end{proposition}

This proposition is the $q=0$ case of a more general $(q,t)$
theorem obtained by Bazlov \cite{bazlov} and Ion \cite{ion}. They consider
$$\tilde{K}(q,t) = \prod_{\alpha \in \rp} \prod_{i \geq 0}
\frac{(1-q^i e^{-\alpha})(1-q^{i+1} e^{\alpha})}{(1-tq^i
  e^{-\alpha})(1-tq^{i+1} e^{\alpha})}$$ 
If $[\tilde{K}(q,t)]$ denotes the constant term
(coefficient of $e^0$) of $\tilde{K}(q,t)$, one defines $\tilde{C}(q,t)
:= \tilde{K}(q,t)/[\tilde{K}(q,t)]$ (upto a minor difference in
convention, this is called the Cherednik kernel in \cite{ion}
). Bazlov and Ion compute the coefficient of
$e^{-\beta}$ in $C(q,t)$ for $\beta$ a positive root of $\kma$.
Their approaches use techniques from Cherednik's theory of Macdonald polynomials.

When $q=0$, $\tilde{C}(0,t)$ reduces to $\xi$ introduced above. 
Though proposition \ref{cl1} is only a special case, it has a
very interesting consequence. 
Ion showed
\cite{ion} that it can be used to give a quick and elegant proof of
the classical fact that for a finite dimensional simple Lie algebra
$\kma$, the partition formed by listing its exponents in descending
order is dual to the partition formed by the numbers of
positive roots at each height (see \S \ref{stwo}). This fact,
first observed empirically by Shapiro and Steinberg was later proved
by Kostant \cite{kost} using his theory of principal three
dimensional subalgebras of $\kma$ and by
Macdonald \cite{macd} via his factorization of the Poincar\'{e} series
of the Weyl group of $\kma$.

The motivation for our approach to proposition \ref{cl1} 
is to thereby obtain a proof of this classical fact via elementary means
(bypassing Macdonald-Cherednik theory).


\section{}\label{stwo}
For completeness sake, we first quickly recall \cite{ion}  how one can
use proposition \ref{cl1} to
deduce the classical fact concerning exponents and heights of roots.

Let $P$ denote the  weight lattice of $\kma$ and $W$ its Weyl group.
For our definition of
exponents, we use the Kostka-Foulkes polynomial $K_{\highroot,0}(t)$ where
$\highroot$ is the highest long root of $\kma$. It is well known that 
this is given by
\beq \label{kost}
K_{\highroot,0}(t) = \sum_{j=1}^l t^{m_j}
\eeq
 where
$m_1,m_2,\cdots,m_l$ are the exponents of $\kma$.
The Kostka-Foulkes polynomials are the elements of the transition matrix
between the Schur and the Hall-Littlewood bases of
$\rationals[t][P]^W$.
 They may be alternatively defined via
Lusztig's $t-$analog of weight multiplicity; we have 
\begin{align}
 K_{\highroot,0}(t) &= \sum_{w \in W} (-1)^{\len(w)}
\, \tko[w(\highroot+\rho) - (0+\rho)]\notag\\
&= \text{coeff. of } e^0 \text{ in  }\frac{\displaystyle\sum_{w \in W} (-1)^{\len(w)} e^{w (\highroot+\rho) -
  \rho}}{\displaystyle\prod_{\alpha \in \rp} (1-te^{-\alpha})}\notag
\end{align}
where $\tko$ is the $t-$analog of Kostant's partition function.

The last equation can be rewritten as :
\begin{align}
K_{\highroot,0}(t) =\text{coeff. of } e^0 \text{ in  } \frac{\displaystyle\sum_{w \in W} (-1)^{\len(w)} e^{w (\highroot+\rho) -
  \rho}}{\displaystyle\prod_{\alpha \in \rp} (1-e^{-\alpha})} \cdot \xi \label{last}\end{align}
where $\xi$ was defined earlier.
The expression in \eqref{last} (from which we need to extract the
coefficient of $e^0$) is just the product $\sch[\highroot]
\, \xi$ where $\sch[\highroot]$ is the formal character of the
adjoint representation of $\kma$. This follows from the Weyl character
  formula and the fact that adjoint representation is irreducible with
  highest weight $\highroot$. Now,
\be
\item  $\sch[\highroot] = le^0 + \displaystyle\sum_{\alpha \in \rp} (e^{\alpha} +
  e^{-\alpha})$ and
\item From equation \eqref{xiq}, the power series for $\xi$ has
  constant term 1 and only
  involves terms of the form $e^{-\gamma}$ for $\gamma \in Q^+$.
\ee
Thus,  $$\text{coeff. of } e^0 \text{ in } \sch[\highroot] \, \xi = l +
  \displaystyle\sum_{\alpha \in \rp} (\text{ coeff. of } e^{-\alpha} \text{ in } \xi)$$
From proposition \ref{cl1}, the right hand side equals $l + \displaystyle\sum_{\alpha \in
  \rp} (t^{\htt(\alpha)} - t^{\htt(\alpha)-1})$. Letting $a_i:=\#\{\beta
  \in \rp: \htt \beta =i\}$, this last sum becomes \\$l + \displaystyle\sum_{i \geq 1}
  a_i (t^i - t^{i-1}) = (a_1-a_2) t + (a_2 - a_3)t^2 + \cdots$ (since
  $a_1=l$). Comparing with equation \eqref{kost}, we get $a_i-
  a_{i+1}$ is the number of times $i$ appears as an exponent of
  $\kma$. This is exactly the classical result. \qed

\section{}
\noindent
{\bf Proof of proposition \ref{cl1}:} Given $\gamma \in Q^+$, let $\p$ be the set of
all partitions of $\gamma$ into a sum of positive roots. Given such a
partition $\pi \in \p$, say $$\pi : \;\;\;\;\;\;\; \gamma = \displaystyle\sum_{\alpha \in
  \rp} c_{\alpha} \alpha \;\;\;\;\;\; \; (c_\alpha \in \integers^{\geq 0})$$
let $n(\pi):=\displaystyle\sum_\alpha c_\alpha$ be the total number of
parts (counting repetitions) and
$d(\pi):=\#\{\alpha: c_\alpha \neq 0\}$ be the number of distinct
parts in $\pi$. From equation \eqref{xiq}, it is clear that
\beq \label{comb} 
\text{Coeff. of } e^{-\gamma} \text{ in } \xi = \displaystyle\sum_{\pi \in \p}
t^{n(\pi)-d(\pi)}(t-1)^{d(\pi)}
\eeq
\underline{Notation:} \be \item Given a subset $A \subset \p$, let $\wt (A) :=
\displaystyle\sum_{\pi \in A} t^{n(\pi)-d(\pi)}(t-1)^{d(\pi)}$. Thus
the coeff. of $e^{-\gamma} \text{ in } \xi$ equals $\wt(\p)$.

\item Given a simple root $\alpha_i$, let 
\begin{align*}
\ppi &:=\{ \pi \in \p:  \alpha_i \text{ occurs as one of the parts in
} \pi\} \\
\ppni &:= \{ \pi \in \p:
  \alpha_i \text{ does not occur as a part in } \pi\}
\end{align*}
\ee

\vspace{2mm}
Let $(\cdot,\cdot)$ denote a nondegenerate, $W-$invariant
  symmetric bilinear form on the dual of the Cartan subalgebra and let
  $s_j \in W$ ($j=1\cdots n$) be the simple reflection correponding to
  $\alpha_j$.

We will sprove proposition \ref{cl1} by induction on $\htt \beta$. If $\htt
  \beta =1$, $\beta$ is a simple root. It is then clear from equation
  \eqref{xiq} that the coefficient of $e^{-\beta}$ is $t-1 = t^1 -
  t^0$. Now suppose $\beta \in \rp$ with $h:=\htt \beta \geq
  2$. Assume the proposition is true for all positive roots of height
  $<h$.  
  Choose a simple root $\alpha_i$ such that
  $(\beta,\alpha_i) >0$ (such $\alpha_i$ exists since
  $(\beta,\beta)>0$). Now $h \geq 2$ implies that $s_i \beta = \beta -
  2\frac{(\beta,\alpha_i)}{(\alpha_i,\alpha_i)} \,\alpha_i$ is 
a positive root of height  $<h$.

\vspace{2mm}
\noindent
\underline{\bf Fact 1:} $\ppni[\beta]$ is in bijection with $\ppni[s_i \beta]$\\
{\em Proof:} Given a partition $\pi: \; \beta = \displaystyle\sum_{\alpha \in \rp}
c_\alpha \alpha \in \ppni[\beta]$, we can form a partition of $s_i
\beta$ as follows:
$$\tilde{\pi}: \;  s_i \beta = \displaystyle\sum_{\alpha \in \rp} c_\alpha \,(s_i
\alpha)$$ Since $\alpha_i$ is not one of the parts of $\pi$, all the
parts of $\tilde{\pi}$ are positive roots, none equal to $\alpha_i$. It is
clear that $\pi \mapsto \tilde{\pi}$ sets up the required bijection. Further, since
$n(\pi) = n(\tilde{\pi})$ and $d(\pi) = d(\tilde{\pi})$, we have
\beq \label{fa1}
\wt (\ppni[\beta]) = \wt(\ppni[s_i \beta])
\eeq

\noindent
\underline{\bf Fact 2:}
\beq \label{two}
\wt(\ppi[\beta]) = t \wt(\ppi[\beta-\alpha_i]) +
(t-1)\wt(\ppni[\beta-\alpha_i])
\eeq
{\em Proof:}  There is an obvious bijection between the sets
$\p[\beta-\alpha_i]$ and $\ppi[\beta]$ obtained by sending a partition
$\pi$ in the first set to the partition $\bar{\pi}$ obtained by
adjoining the extra part $\alpha_i$ to $\pi$. In order to see how $\wt(\pi)$
compares with $\wt(\bar{\pi})$, we write $\p[\beta-\alpha_i] =
\ppi[\beta-\alpha_i] \cup \ppni[\beta-\alpha_i]$. For $\pi \in
\ppi[\beta-\alpha_i]$, the extra part $\alpha_i$ in $\bar{\pi}$ is a
repeat part and thus $$\wt(\bar{\pi}) = t \, \wt(\pi)$$ while for
 $\pi \in \ppni[\beta-\alpha_i]$, the extra $\alpha_i$ in $\bar{\pi}$
 is a new distinct part and thus $\wt(\bar{\pi}) = (t-1) \,
 \wt(\pi)$. This proves equation \eqref{two}. \qed

Let $k:=2\frac{(\beta,\alpha_i)}{(\alpha_i,\alpha_i)} >0$ and consider the
$\alpha_i-$string through $\beta$: $$\beta,\;\beta - \alpha_i, \cdots,
\;\beta - k \alpha_i$$ 
Each of these is a positive root.
We now rewrite equation \eqref{two} as 
$$\wt(\ppi[\beta]) - \wt(\ppi[\beta-\alpha_i]) = (t-1) \wt(\p[\beta-\alpha_i])  $$
Iterating this equation  $k$ times with $\beta - j
\alpha_i$ in place of $\beta$ ($0 \leq j \leq k-1$) and summing the
resulting equations, we get
$$\wt(\ppi[\beta])-  \wt(\ppi[\beta-k\alpha_i]) = (t-1) \displaystyle\sum_{j=1}^k  \wt(\p[\beta-j\alpha_i]) $$

By induction hypothesis, $$\wt(\p[\beta-j\alpha_i]) = t^{h-j} -
t^{h-j-1}$$
Further, 
\begin{align*}
\wt(\ppi[\beta-k\alpha_i]) &=  \wt(\p[\beta-k\alpha_i]) -
\wt(\ppni[\beta-k\alpha_i])\\
&=(t^{h-k} -t^{h-k-1}) - \wt(\ppni[\beta-k\alpha_i])
\end{align*}
Since $\beta - k \alpha_i = s_i \beta$, we can use equation
\eqref{fa1}. This  gives
\begin{align*}
\wt(\ppi[\beta]) + \wt(\ppni[\beta]) &=  (t-1) \displaystyle\sum_{j=1}^k ( t^{h-j}
-t^{h-j-1}) + (t^{h-k} -t^{h-k-1}) \\
&= t^h - t^{h-1}
\end{align*}
Since the left hand side equals $\wt(\p[\beta])$, proposition \ref{cl1} is
proved. \qed

\noindent
{\em Acknowledgements} : The author would like to thank John Stembridge
 for   bringing references \cite{bazlov} and \cite{ion} to his
 attention and  for his comments on an earlier draft of this note.